\documentclass[12pt]{article} 
\usepackage{amssymb,amsthm,amsmath}
\usepackage{caption}
\usepackage{subcaption}
\usepackage{colortbl}
\usepackage{xcolor}

\newcommand{\comment}[1]{}

\newcommand{\entryi}{\alpha}
\newcommand{\entryii}{\beta}
\newcommand{\entryiii}{\gamma}
\newcommand{\ttt}[1]{{\em #1}}

\newtheorem{theorem}{Theorem}[section]

\title{Tic-Tac-Toe on an Affine Plane of order 4}

\author{
Peter Danziger\thanks{Department of Mathematics, Ryerson University, Toronto, Ontario, Canada}
\and
Melissa A. Huggan\thanks{Department of Mathematics and Computer Science, Mount Allison University, Sackville, New Brunswick, Canada}
\and
Rehan Malik$^*$
\and
Trent G.\ Marbach$^*$
}

\begin{document}
\maketitle

\begin{abstract}
The game of \ttt{Tic-Tac-Toe} is well known. In particular, in its classic version it is famous for neither player having a winning strategy.
While classically it is played on a grid, it is natural to consider the effect of playing the game on richer structures, such as finite planes.
Playing the game of \ttt{Tic-Tac-Toe} on finite affine and projective planes has been  previously studied.
While the second player can usually force a draw, for small orders the first player has a winning strategy.
In this regard, a computer proof that \ttt{Tic-Tac-Toe} played on the affine plane of order 4 is a first player win has been claimed.  
In this note we use techniques from the theory of latin squares and transversal designs to give a human verifiable, explicit proof of this fact. 
\end{abstract}

\comment{
}
\section{Introduction}

The game of \ttt{Tic-Tac-Toe} is well known, particularly for there being no winning strategy.
While classically it is played on a grid, it is natural to consider the effect of playing the game on richer structures, such as finite planes.
Carroll and Dougherty \cite{CarrollDougherty} examined \ttt{Tic-Tac-Toe} played on finite affine planes (denoted $\pi_n$)  and finite projective planes (denoted $\Pi_n$). 
They showed that for $n\leq 4$, $\pi_n$ is a first player win, whereas for $n > 4$, it is always possible for the second player to force a draw.
In this regard, they claim a computer proof that $\pi_4$ is a first player win. They note that despite this, 
they are able to hold a local \ttt{Tic-Tac-Toe} tournament played on $\pi_4$ as the strategy is not obvious.
In this paper we give an explicit, 
human verifiable proof that $\pi_4$ is a first player win.
We do this by using the language of latin squares and transversal designs, see \cite{Handbook,Dembowski} for further definitions and details. 

The game of \ttt{Tic-Tac-Toe} is a positional game. 
Positional games are games where players alternately pick points and a player wins if they are the first to occupy all points of some specified configuration(s), otherwise the game is a draw. 
Positional games, including classical \ttt{Tic-Tac-Toe} played on a square grid, have been well studied (see~\cite{Beck05,Beck08, ES73}). 
In addition, playing the game on other structures has been investigated, for example on hypercubes \cite{GH02}, graphs \cite{Beeler18}, 
as well as affine and projective planes \cite{CarrollDougherty}.
Other variants, including quantum moves~\cite{Goff2006}, infinite boards~\cite{CCIKLST15}, and movable pieces~\cite{LR17}, have also been explored.

In the context of playing \ttt{Tic-Tac-Toe} on a finite affine plane, players alternate choosing points on the plane; the first player who collects all of the points of a line wins the game. Following the convention of \cite{CarrollDougherty}, we call the first player Xeno and the second player Ophelia. 

Classically, an affine plane, $\pi_n$, is a set of $n^2$ {\em points} and $n(n+1)$ {\em lines}, each of which contains an $n$-set of points, with the property that given any pair of points there is exactly one line that contains them both. Further, a set of lines that partition the point set is called a {\em parallel class} and the lines of $\pi_n$ can be partitioned into parallel classes. It is well known that $\pi_n$ exists for every prime power order, see~\cite{Handbook}, but it is a 
longstanding open problem whether $\pi_n$ exists for any non-prime order $n$.
In this paper we focus on the case when $n=4$; so $\pi_4$ consists of 16 points and 20 lines, each line contains 4 points, 
and the lines can be divided into 5 parallel classes. In order to introduce $\pi_4$ more explicitly we use the language of latin squares and transversal designs, which we describe below.

A {\em latin square} of order $n$ is an $n\times n$ array filled with elements from an $n$-set, called the symbol set, such that each symbol appears exactly once in each row and once in each column. Two latin squares of order $n$ are called {\em orthogonal} if every pair of symbols appears in the superimposition of the two squares. Given $k$ latin squares, they are {\em mutually orthogonal} if they are pairwise orthogonal and we refer to them as $k$-MOLS$(n)$.
Figure~\ref{fig:2MOLS4} gives three mutually orthogonal latin squares of order 4 ($3$-MOLS$(4)$).
We refer to the $i^{\rm th}$ row as $r_i$, the $i^{\rm th}$ column as $c_i$, the $i^{\rm th}$ symbol in the first square as $\entryi_i$, the $i^{\rm th}$ symbol of the second square as $\entryii_i$, and the $i^{\rm th}$ symbol of the third square as $\entryiii_i$. 
For example, in Figure~\ref{fig:2MOLS4} the symbol in $r_3, c_2$ is  $\entryi_4$ in the first square, $\entryii_1$ in the second square, and $\entryiii_3$ in the third square. 
We refer to the corresponding 5-tuple as an {\em entry} of the three mutually orthogonal squares. For example, the entry in row $r_3$, column $c_2$ of Figure~\ref{fig:2MOLS4} is the 5-tuple $(r_3, c_2, \entryi_4, \entryii_1,\entryiii_3)$. 

\begin{figure}[h]
\[
\begin{array}{ccccc}
         \begin{array}{|c|c|c|c|c|c|}
         \hline
         \entryi_{1}&\entryi_{2}&\entryi_{3}&\entryi_{4}\\
         \hline
         \entryi_{2}&\entryi_{1}&\entryi_{4}&\entryi_{3}\\
         \hline
          \entryi_{3}&\entryi_{4}&\entryi_{1}&\entryi_{2}\\
         \hline
          \entryi_{4}&\entryi_{3}&\entryi_{2}&\entryi_{1}\\
         \hline
         \end{array}
& \hspace*{2em} &
        \begin{array}{|c|c|c|c|c|c|}
         \hline
         \entryii_1& \entryii_2& \entryii_3& \entryii_4\\
         \hline
         \entryii_4& \entryii_3& \entryii_2& \entryii_1\\
         \hline
         \entryii_2& \entryii_1& \entryii_4& \entryii_3\\
         \hline
         \entryii_3& \entryii_4& \entryii_1& \entryii_2\\
         \hline
         \end{array}
& \hspace*{2em} &
          \begin{array}{|c|c|c|c|c|c|}
         \hline
         \entryiii_1& \entryiii_2& \entryiii_3& \entryiii_4\\
         \hline
         \entryiii_3& \entryiii_4& \entryiii_1& \entryiii_2\\
         \hline
         \entryiii_4& \entryiii_3& \entryiii_2& \entryiii_1\\
         \hline
         \entryiii_2& \entryiii_1& \entryiii_4& \entryiii_3\\
         \hline
         \end{array}  
\end{array}
\]
        \caption{Three mutually orthogonal latin squares of order 4.}
        \label{fig:2MOLS4}
\end{figure}

Given integers $k$ and $n$, a {\em Transversal Design}, TD$(k,n)$, is a triple $(X, {\cal G}, {\cal B})$, where $X$ is a $kn$-set of {\em points}, ${\cal G}$ is a set of $n$-sets of $X$ that partition $X$, called {\em groups} and ${\cal B}$ is a collection of $k$-sets of $X$, called {\em blocks}, such that every pair of points either appears in a block $B\in {\cal B}$, or in $G\in{\cal G}$, but not both.
It is well known that $k$-MOLS$(n)$ is equivalent to a TD$(k+2,n)$.
For example, taking the entries (5-tuples) corresponding to all of the 16 cells of the $3$-MOLS$(4)$ in Figure~\ref{fig:2MOLS4} above as blocks gives us the blocks of a transversal design, TD$(5,4)$. 

A {\em parallel class} of a transversal design is a disjoint set of blocks whose union contains all of the points. If the blocks of a transversal design can be partitioned into parallel classes, it is called {\em resolvable} and we refer to a {\em Resolvable Transversal Design}, RTD$(k,n)$. Given $k$-MOLS($n$) it is well known that we can use the entries of the final square to index the entries in each parallel class to obtain a resolvable transversal design, RTD$(k+1,n)$.

The RTD(4,4) corresponding to the three latin squares from Figure~\ref{fig:2MOLS4} has 16 points, 16 blocks and 4 groups.
The point set, $X$, consists of the 16 points defined by the rows, columns, and the symbol sets of the first two squares

$$X = \{r_i, c_j, \entryi_k, \entryii_\ell : 1\leq i, j, k, \ell\leq 4\}.$$ 
The blocks are defined by the 4-tuple corresponding to each cell from the first two squares,
see Figure~\ref{fig: TD(4,4)}. 
Each of the rows in Figure~\ref{fig: TD(4,4)} is a parallel class corresponding to the cells containing $\entryiii_i$.
The sets of rows, columns, entries of the first square, and entries of the second square are the groups of the transversal design, see Figure~\ref{groups}.

\begin{figure}[h]
\[
\begin{array}{lcccc}
\entryiii_1: & \{r_1, c_1, \entryi_1, \entryii_1\}, & \{r_4, c_2, \entryi_3, \entryii_4\}, & \{r_2, c_3, \entryi_4, \entryii_2\}, & \{r_3, c_4, \entryi_2, \entryii_3\}; \\
\entryiii_2: &
\{r_3, c_3, \entryi_1, \entryii_4\}, & \{r_1, c_2, \entryi_2, \entryii_2\}, & 
\{r_4, c_1, \entryi_4, \entryii_3\}, & \{r_2, c_4, \entryi_3, \entryii_1\}; \\
\entryiii_3: &
\{r_2, c_1, \entryi_2, \entryii_4\}, & \{r_3, c_2, \entryi_4, \entryii_1\}, & 
\{r_1, c_3, \entryi_3, \entryii_3\}, &\{r_4, c_4, \entryi_1, \entryii_2\}; \\
\entryiii_4: &
\{r_3, c_1, \entryi_3, \entryii_2\}, & \{r_2, c_2, \entryi_1, \entryii_3\}, & \{r_4, c_3, \entryi_2, \entryii_1\}, & \{r_1, c_4, \entryi_4, \entryii_4\}. \\
\end{array}
\]
\caption{The blocks of the resolvable transversal design corresponding to the three latin squares in Figure~\ref{fig:2MOLS4}.}
\label{fig: TD(4,4)}
\end{figure}

We may obtain $\pi_4$ by taking the point set to be $X$ above, and as lines we take the blocks of the RTD(4,4) in Figure~\ref{fig: TD(4,4)} along with the groups from Figure~\ref{groups}, which form an additional parallel class.
We refer to this additional parallel class as the {\em index} parallel class and refer to the blocks of this class as the {\em row block}, {\em column block}, {\em symbol set one block}, and {\em symbol set two block} respectively (see Figure~\ref{groups}).
\begin{figure}[h]
\[
\begin{array}{cccc}
\{r_1, r_2, r_3, r_4\}, & \{c_1, c_2, c_3, c_4\}, & \{\entryi_1, \entryi_2, \entryi_3, \entryi_4\}, & \{\entryii_1, \entryii_2, \entryii_3, \entryii_4\}. 
\end{array}
\]
\caption{\label{groups}
The index class of $\pi_4$, equivalently the groups of the RTD(4,4).}
\end{figure}

In order to play \ttt{Tic-Tac-Toe} on $\pi_4$, players move by alternately choosing points, and win if they complete a line.
We can view the game as being played on a pair of MOLS(4), picking rows, columns, and symbols. A player wins if they have chosen all of the 
components of an entry or chosen all of the rows, all of the columns, or all of the symbols from a square.

We note that when playing on $\pi_4$, the third square is suppressed and we are only playing on the first two squares.
Indeed, playing on the TD$(5,4)$ together with the index class, equivalently playing on all three MOLS$(4)$,
is the same as playing on the projective geometry $\Pi_4$, for which it is known that Ophelia can force a draw~\cite{CarrollDougherty}.

We can represent this in a more standard \ttt{Tic-Tac-Toe} grid-like structure as follows. We create a square grid in which the blocks of $\gamma_1$ are the rows and the blocks of $\gamma_2$ are the columns as shown in Figure~\ref{grid}. 
We define a diagonal to be a set of four cells such that each row and column is represented exactly once.
Now, choosing any of the following is a winning set of cells:
\begin{itemize}
\item
a row (corresponding to a block from $\gamma_1$), 
\item
a column (corresponding to a block from $\gamma_2$), 
\item
a diagonal containing different index types, see Figure~\ref{diagonals} (corresponding to a block of $\gamma_3$ or $\gamma_4$),
\item
all cells containing the same index type (corresponding to a block from the index class).
\end{itemize}

\renewcommand{\arraystretch}{1.5}
\begin{figure}
\[
\begin{array}{|c|c|c|c|} \hline
r_2 & \entryi_4 & \entryii_2 & c_3 \\ \hline
\entryii_1 & c_1 & r_1 &  \entryi_1 \\ \hline
c_4 &  \entryii_3 &  \entryi_2 & r_3 \\ \hline
\entryi_3 & r_4 & c_2 &  \entryii_4 \\ \hline
\end{array}
\]
\caption{\label{grid} A \ttt{Tic-Tac-Toe} grid.}
\end{figure}
\renewcommand{\arraystretch}{1}
\begin{figure}
\[
\begin{array}{ccc}
\begin{array}{|c|c|c|c|} \hline
x & y & z & w \\ \hline
y & x & w & z \\ \hline
z & w &  x & y \\ \hline
w & z & y & x  \\ \hline
\end{array}
& & 
\begin{array}{|c|c|c|c|} \hline
x & y & z & w \\ \hline
w & z & y & x \\ \hline
y & x & w & z \\ \hline
z & w & x & y  \\ \hline
\end{array}
\\
\\
(a) \mbox{ The winning diagonals from } \gamma_3.
& & 
(b) \mbox{ The winning diagonals from } \gamma_4.
\end{array}
\]
\caption{
\label{diagonals}
The winning diagonals corresponding to $\gamma_3$ and $\gamma_4$.}
\end{figure}

A {\em paratopism} of 2-MOLS($n$) is a map, which can be described by a 5-vector $(\pi,\sigma_r,\sigma_c, \sigma_{\entryi}, \sigma_{\entryii})$, where $\pi$ is a permutation that maps rows, columns, and symbols in the squares between themselves and $\sigma_r,\sigma_c, \sigma_{\entryi}, \sigma_{\entryii}$ are permutations of the resulting rows, columns, symbols in the first square, and symbols in the second square respectively.
We use $i$ to represent the identity permutation.
An {\em autoparatopism} is a paratopism that leaves the squares unchanged, see \cite{EW16}.
Clearly, any autoparatopism of a pair of MOLS($4$) generates an isomorphism of $\pi_4$.
It is well known that $\pi_4$ is unique up to isomorphism, see \cite{Handbook, EW16}.
Thus the squares in Figure~\ref{fig:2MOLS4} and the corresponding RTD in Figure~\ref{fig: TD(4,4)} are  also unique up to isomorphism.

We can describe a game of \ttt{Tic-Tac-Toe} on $\pi_4$ by giving the sequence of moves in order. For ease of reading we give Ophelia's moves in parentheses. Thus a game on $\pi_4$ might look like:
$$r_1, (r_2), r_3, (c_1), \entryi_2, (r_4), c_2, ({\entryii_2}), \entryi_4, ({\entryii_1}), c_4, (\entryii_3), \entryii_4;$$ 
which results in a Xeno win as he now has the line $(r_1, c_4, \entryi_4, \entryii_4)$.
Any move whose omission will result in the opponent winning the game
is called a {\em forced move} and we indicate this type of move by placing a line over it.
In the above example $\entryii_2$ is forced because if Ophelia does not take $\entryii_2$ at this point, Xeno will take it and win with the line $(r_1, c_2, \entryi_2, \entryii_2)$.

Note that in the example above, if Ophelia does not make the move $\entryii_3$ and instead takes $\entryii_4$ in an attempt to stop Xeno from winning, Xeno can then take $\entryii_3$ and he will still win the game with the line $(r_3,c_4,\entryi_2,\entryii_3)$. 
Generalizing this notion, if it is Ophelia's move and she has two forced moves $x$ and $y$ (as with $\entryii_3$ and $\entryii_4$ above), Xeno is able to win the game regardless of Ophelia's choice and so the game ends with a Xeno win. We denote this situation by $X_W(x,y)$. Thus the sequence of moves in the game above would be written as
$$r_1, (r_2), r_3, (c_1), \entryi_2, (r_4), c_2, (\overline{\entryii_2}), \entryi_4, (\overline{\entryii_1}), c_4, X_W(\entryii_3, \entryii_4).$$

Since we wish to describe a winning strategy for Xeno, we only need to give optimal moves for Xeno and thus we will not always indicate that a Xeno move is forced for ease of exposition. On the other hand, if Ophelia's move is not forced, we either need to exploit a symmetry of the game or enumerate strategies for all possible responses she could make, or a combination of the two. 

At any point in the game we can describe the game thus far by listing the points that Xeno has chosen and the points that Ophelia has chosen thus far. If Xeno has chosen the set of points $A$ and Ophelia has chosen the set of points $B$, we denote these by $X =A$ and $O=B$ respectively. 

\section{The Affine Plane $\pi_4$ is a Xeno win}

We are now ready to present our main Theorem.
\begin{theorem} 
The Affine Plane of order $4$, $\pi_4$, is a first player (Xeno) win.
\end{theorem}

\begin{proof}
We start by playing on an unlabelled $\pi_4$. As the game progresses we label the points so that they are consistent with the first two latin squares in Figure~\ref{fig:2MOLS4} and thus the RTD in Figure~\ref{fig: TD(4,4)} together with the index class in Figure~\ref{groups}.

Xeno initially plays on a point which we arbitrarily label $r_1$. Ophelia now plays on another point, which we label $r_2$.  Now, the line containing $r_1$ and $r_2$ is the 
row block of the corresponding RTD. This block is contained in a parallel class $P$, which we take to be the index class and so the blocks of that parallel class will label the columns, symbol ones ($\entryi_i\,$) and symbol twos ($\entryii_i$) in some order. Xeno now plays on another point of the row block, which we take to be $r_3$.

If Ophelia's next move is in the row block, it must be $r_4$. In this case, Xeno now plays a point in another block, which we label $c_1$ and thus the block of $P$ that contains it is the column block.
Alternatively, if Ophelia's next move is not in the row block (i.e.\ not $r_4$), we label her move $c_1$ and thus the block of $P$ that contains it is the column block.
Xeno now plays a point on the block through $r_2$, $c_1$, which we label $\entryi_2$ to be consistent with the latin squares and RTD given in Figures~\ref{fig:2MOLS4} and \ref{fig: TD(4,4)} respectively.
The fourth point in this block is thus labelled $\entryii_4$. 
Note that this argument implicitly relies on the fact that there is only one $\pi_4$ up to isomorphism.
Thus, the configuration to this point is either 
$$X= (r_1, r_3, c_1),\; O = (r_2, r_4) \mbox{ or } X= (r_1, r_3, \entryi_2),\; O = (r_2, c_1).$$
We consider each of these possibilities in turn.

\begin{figure}[h]
     \centering
     \begin{subfigure}[b]{0.3\textwidth}
         \centering
         \begin{tabular}{|c|c|c|c|c|c|}
         \hline
         \cellcolor{red}{1} & \cellcolor{green}{2} & \cellcolor{green}{3} &  \cellcolor{green}{4}\\
         \hline
         \cellcolor{cyan}{2} & \cellcolor{yellow}{1} & \cellcolor{yellow}{4} & \cellcolor{yellow}{3}\\
         \hline
         \cellcolor{red}{3} & \cellcolor{green}{4} & \cellcolor{green}{1} & \cellcolor{green}{2}\\
         \hline
         \cellcolor{cyan}{4} & \cellcolor{yellow}{3} & \cellcolor{yellow}{2} & \cellcolor{yellow}{1}\\
         \hline
         \end{tabular}
     \end{subfigure}
     \begin{subfigure}[b]{0.3\textwidth}
         \centering
                 \begin{tabular}{|c|c|c|c|c|c|}
         \hline
          \cellcolor{red}{1} &  \cellcolor{green}{2} & \cellcolor{green}{3} & \cellcolor{green}{4}\\
         \hline
        \cellcolor{cyan}{4} & \cellcolor{yellow}{3} & \cellcolor{yellow}{2} &  \cellcolor{yellow}{1}\\
         \hline
          \cellcolor{red}{2} & \cellcolor{green}{1} & \cellcolor{green}{4} & \cellcolor{green}{3}\\
         \hline
        \cellcolor{cyan}{3} & \cellcolor{yellow}{4} & \cellcolor{yellow}{1} & \cellcolor{yellow}{2}\\
         \hline
         \end{tabular}
     \end{subfigure}
        \caption{The 2-MOLS(4) in the first case. Entries in the first square correspond to $\entryi_i$ and the second to the $\entryii_i$.
        Cells are coloured: Blue - Neither can win with that entry; Green - Xeno controls one coordinate; Red - Xeno controls two coordinates; Yellow - Ophelia controls one coordinate.}
        \label{fig:2MOLS4coloured}
\end{figure}

The first case is illustrated in Figure~\ref{fig:2MOLS4coloured}.
In this case we find two autoparatopisms $(\pi,\sigma_r,\sigma_c, \sigma_{\entryi}, \sigma_{\entryii})$ that fix the moves of Xeno and Ophelia, given by:
\[
\psi_1 = (i,\, (13)(24),\, i,\, (13)(24),\, (12)(34)),
\]
\mbox{ and }
\[
\psi_2 = ((34),\,(13),\, (23),\, (2134),\,(3124)).
\]

Referring to Figure~\ref{fig:2MOLS4coloured}, any autoparatopism must map the same colors to each other.
This means that we cannot interchange rows or columns with each other or with symbols and so the only option is to either swap the two symbol sets (squares) or to leave them alone, in which case we then derive $\psi_2$ and $\psi_1$ respectively.

We use these autoparatopisms to map between different games to show that seemingly different games are actually the same game.

Up to these autoparatopisms Ophelia's only unique moves are on $c_2, c_4,\entryi_1$, and $\entryi_2$. If Ophelia now plays $c_2$ or $c_4$ then a winning sequence for Xeno is as follows: 
$$\entryi_1, (\overline{\entryii_1}), \entryi_3, (\overline{\entryii_2}), c_3, X_W(\entryii_3, \entryii_4).$$ 
If Ophelia plays an autoparatopism $\psi_i$, $i=1,2$, of $c_2$ or $c_4$, then play proceeds by the autoparatopic sequence of moves:
\[\psi_i(\entryi_1), (\overline{\psi_i(\entryii_1})), \psi_i(\entryi_3), (\overline{\psi_i(\entryii_2})), \psi_i(c_3), X_W(\psi_i(\entryii_3), \psi_i(\entryii_4)).\]

If Ophelia plays $\entryi_2$, then a winning sequence for Xeno is as follows:  
$$\entryi_1, (\overline{\entryii_1}), \overline{c_3}, (\overline{\entryii_4}), \entryi_3, X_W(\entryii_2, \entryii_3).$$
Similarly to the case above, if Ophelia plays an autoparatopism $\psi_i(\entryi_2)$, $i=1,2$, then play proceeds by the autoparatopic sequence of moves:
\[\psi_i(\entryi_1), (\overline{\psi_i(\entryii_1})), \psi_i(c_3), (\overline{\psi_i(\entryii_4})), \psi_i(\entryi_3), X_W(\psi_i(\entryii_2), \psi_i(\entryii_3)).\]

If Ophelia plays $\entryi_1$, play proceeds as $\entryii_2$, ($\overline{\entryi_3}$), $\entryii_4$.
The remaining moves for Ophelia and a winning response for Xeno in each case are given in 
Table~\ref{Moves-Case 1}.
If Ophelia plays an autoparatopism, $\psi_i$, $i=1,2$, applied to $\entryi_1$, play proceeds by applying $\psi_i$ to the appropriate sequence of moves. 
\pagebreak

\addtocounter{table}{6}
\addtocounter{figure}{1}
\begin{table}[h]
\[
\begin{array}{|c|c|c|c|c|} \hline
O & X & O & X & X_W\\ \hline
\vspace*{-1ex} & & & & \\
c_2 & \overline{\entryii_3} & \overline{\entryii_1} & \overline{c_4} & X_W(\entryi_2, \entryi_4) \\[1ex] 
c_3 & c_4 & \overline{\entryi_4} & \overline{\entryi_2} & X_W(c_2, \entryii_3) \\[1ex] 
c_4 & \overline{\entryii_1} & \overline{\entryii_3} & \overline{c_2} & X_W(\entryi_2, \entryi_4) \\[1ex] 
\entryi_2 & \overline{\entryi_4} & \overline{c_4} & \overline{\entryii_1} & X_W(c_2, \entryii_3) \\[1ex] 
\entryi_4 & \overline{\entryi_2} & \overline{c_2} & \overline{\entryii_3} & X_W(c_4, \entryii_1) \\[1ex] 
\entryii_1 & \overline{c_4} & \overline{\entryi_4} & \overline{\entryi_2} & X_W(c_2, \entryii_3) \\[1ex] 
\entryii_3 & \overline{c_2} & \overline{\entryi_2} & \overline{\entryi_4} & X_W(c_4, \entryii_1) \\[1ex] 
\hline
\end{array}
\]
\caption{\label{Moves-Case 1} Case 1: Remaining Ophelia moves and Xeno's responses.}
\end{table}

\begin{figure}[h]
     \centering
     \begin{subfigure}[b]{0.3\textwidth}
         \centering
         \begin{tabular}{|c|c|c|c|c|c|}
         \hline
         \cellcolor{cyan}{1} & \cellcolor{red}{2} & \cellcolor{green}{3} &  \cellcolor{green}{4}\\
         \hline
         \cellcolor{cyan}{2} & \cellcolor{yellow}{1} & \cellcolor{yellow}{4} & \cellcolor{yellow}{3}\\
         \hline
         \cellcolor{cyan}{3} & \cellcolor{green}{4} & \cellcolor{green}{1} & \cellcolor{red}{2}\\
         \hline
         \cellcolor{yellow}{4} & {3} & \cellcolor{green}{2} & {1}\\
         \hline
         \end{tabular}
     \end{subfigure}
     \begin{subfigure}[b]{0.3\textwidth}
         \centering
                 \begin{tabular}{|c|c|c|c|c|c|}
         \hline
          \cellcolor{cyan}{1} &  \cellcolor{red}{2} & \cellcolor{green}{3} & \cellcolor{green}{4}\\
         \hline
        \cellcolor{cyan}{4} & \cellcolor{yellow}{3} & \cellcolor{yellow}{2} &  \cellcolor{yellow}{1}\\
         \hline
          \cellcolor{cyan}{2} & \cellcolor{green}{1} & \cellcolor{green}{4} & \cellcolor{red}{3}\\
         \hline
        \cellcolor{yellow}{3} & {4} & \cellcolor{green}{1} & {2}\\
         \hline
         \end{tabular}
     \end{subfigure}
        \caption{The 2-MOLS(4) in the second case. Entries in the first square correspond to $\entryi_i$ and the second to the $\entryii_i$.
        Cells are coloured: Blue - Neither can win with that entry; Green - Xeno controls one coordinate; Red - Xeno controls two coordinates; Yellow - Ophelia controls one coordinate.}
        \label{fig:2MOLS4coloured2}
\end{figure}

We now consider the case where $X= (r_1, r_3, \entryi_2),\; O = (r_2, c_1)$, which is illustrated in Figure~\ref{fig:2MOLS4coloured2}.
As before, any autoparatopism must map the same colours to each other.
This means that we cannot interchange rows or columns with each other or with symbols and so the only option is to either swap the two symbol sets (squares) or to leave them alone. 
It is not hard to see that in either case there is no mapping of the rows, columns and symbols between themselves that preserves both the squares and the colouring.
Thus there are no autoparatopisms fixing Xeno and Ophelia's moves.
However, if Ophelia plays $r_4, c_3, \entryi_1$ or $\entryi_3$, then a winning play for Xeno is as follows:
$$c_2, (\overline{\entryii_2}), \entryi_4, (\overline{\entryii_1}), c_4, X_W(\entryii_3, \entryii_4).$$

The remaining moves for Ophelia and a winning response for Xeno in each case are given in Table~\ref{Moves-Case 2}.

\addtocounter{table}{1}
\begin{table}[h]
\renewcommand{\arraystretch}{1.5}
\[
\begin{array}{|c|c|c|c|c|c|c|c|c|} \hline
O & X & O & X & O & X & O & X & X_W \\ \hline
c_2 & \entryii_3 & \overline{c_4} & \overline{c_3} & \overline{\entryi_3} & \overline{\entryii_1} & \overline{r_4} & \overline{\entryii_4} & X_W(\entryi_1,\entryii_2) \\ \hline  
c_4 & c_2 & \overline{\entryii_2} & c_3 & r_4 & \overline{\entryi_1} & \overline{\entryii_4} & \entryi_3 & X_W(\entryi_4, \entryii_3) \\ 
&&&           & \entryi_1, \entryi_4 & r_4 & \overline{\entryii_1} & \overline{\entryi_3} & X_W(\entryii_3, \entryii_4) \\ 
&&&            & \entryi_3, \entryii_3, \entryii_4 & \entryii_1 & & & X_W(r_4,\entryi_4) \\
&&&            & \entryii_1 & \overline{\entryi_3} & \overline{\entryii_3}& \overline{\entryii_4 } & X_W(r_4, \entryi_1) \\
 \hline
\entryi_4 & \entryii_2 & \overline{c_2} & \entryii_3 & \overline{c_4} & \overline{c_3} & \overline{\entryi_3} & \overline{\entryii_1} & X_W(r_4,\entryii_4) \\ \hline
\entryii_1 & c_4 & \overline{\entryii_3} & \entryi_4 & \overline{\entryii_4} & \overline{\entryii_2} & \overline{c_2} & \overline{\entryi_1} & X_W(r_4, \entryi_3)
\\ \hline
\entryii_2 & c_3 & 
\begin{array}{c}
r_{4}, c_2, \entryi_{1}, \\[-1.3ex]
\entryi_{4}, \entryii_{1}, \entryii_4
\end{array}
& \entryii_3& &&&&X_{W}(c_{4}, \entryi_{3})   \\
&&c_{4}&\entryi_{3}&\overline{\entryii_{3}}&\entryi_{1}&&&X_{W}(\entryi_{4}, \entryii_{4})\\ 
&&\entryi_{3}&\entryii_{4}&\overline{\entryi_{1}}&c_{4}&&&X_{W}(\entryi_{4}, \entryii_{3})\\ 
&&\entryii_{3}&\entryi_{1}&\overline{\entryii_{4}}&\overline{\entryii_{1}}&\overline{r_{4}}&\overline{\entryi_{4}}&X_{W}(c_{2}, \entryi_{3})\\ \hline
\entryii_3 & c_2 & \overline{\entryii_{2}}&\entryi_{4}&\overline{\entryii_{1}}&\overline{\entryii_{4}}&\overline{c_{4}}&\overline{\entryi_{3}}&X_{W}(r_{4},\entryi_{1})\\ \hline
\entryii_4 & \entryii_3 & \overline{c_{4}} & c_{2} & \overline{\entryii_{2}} & \entryi_{4} & \overline{\entryii_{1}} & \overline{\entryi_{3}} & X_{W}(c_{3},\entryi_{1})
 \\ \hline
\end{array}
\]
\caption{\label{Moves-Case 2} Case 2: Remaining Ophelia moves and Xeno's responses.}
\renewcommand{\arraystretch}{1}
\end{table}
\end{proof}

We note that despite the proof we have given here, the tournaments described in \cite{CarrollDougherty} are probably safe for the time being. As well as memorizing the tables above (certainly possible for a dedicated player), applying this strategy would require implementing the autoparatopic responses required, including the mappings implicit in the initial labeling described in the beginning of the proof. This becomes even more difficult if the game is played on a grid as in Figure~\ref{grid}. In this case, the player would also have to translate between the grid and the corresponding Transversal Design.

The techniques applied here involve a detailed analysis of the play. It is unlikely that this will be possible for structures that are much larger than those considered here, though related structures of a similar size might be amenable to such an analysis.
It is natural to ask what the result of playing \ttt{Tic-Tac-Toe}
on other transversal designs might be. We have ascertained by applying similar methods that when playing \ttt{Tic-Tac-Toe} on a TD(4,4) ($\pi_4$ with the index class removed), Ophelia can force a draw. This means that the removal of the four index blocks changes the outcome of the game. 
Interestingly, removing only one index block is still a Xeno win, but removing two or more means that Ophelia can force a draw.

\section*{Acknowledgments}

The first author is supported by the Natural Sciences and Engineering Research Council of Canada, NSERC Discovery Grant number RGPIN-2016-04178. While the research took place the second author was supported by an NSERC Postdoctoral Fellowship PDF-532564-2019 and the fourth author was supported by a Postdoctoral Fellowship from the Fields Institute for Research in Mathematical Sciences; both hosted by the Department of Mathematics at Ryerson University.

\end{document}